\documentclass{amsart}

\usepackage{amsfonts,amssymb,verbatim,amsmath,amsthm,latexsym,textcomp,amscd}
\usepackage{latexsym,amsfonts,amssymb,epsfig,verbatim}
\usepackage{amsmath,amsthm,amssymb,latexsym,graphics,textcomp}
\usepackage{graphicx}
\usepackage{color}
\usepackage{url}

\input xy
\xyoption{all}

\begin{document}

\newtheorem{theorem}{Theorem}[section]
\newtheorem{prop}[theorem]{Proposition}
\newtheorem{lemma}[theorem]{Lemma}
\newtheorem{cor}[theorem]{Corollary}
\newtheorem{definition}[theorem]{Definition}
\newtheorem{conj}[theorem]{Conjecture}
\newtheorem{rmk}[theorem]{Remark}
\newtheorem{claim}[theorem]{Claim}
\newtheorem{defth}[theorem]{Definition-Theorem}

\newcommand{\boundary}{\partial}
\newcommand{\C}{{\mathbb C}}
\newcommand{\integers}{{\mathbb Z}}
\newcommand{\natls}{{\mathbb N}}
\newcommand{\ratls}{{\mathbb Q}}
\newcommand{\bbR}{{\mathbb R}}
\newcommand{\proj}{{\mathbb P}}
\newcommand{\lhp}{{\mathbb L}}
\newcommand{\tube}{{\mathbb T}}
\newcommand{\cusp}{{\mathbb P}}
\newcommand\AAA{{\mathcal A}}
\newcommand\BB{{\mathcal B}}
\newcommand\CC{{\mathcal C}}
\newcommand\DD{{\mathcal D}}
\newcommand\EE{{\mathcal E}}
\newcommand\FF{{\mathcal F}}
\newcommand\GG{{\mathcal G}}
\newcommand\HH{{\mathcal H}}
\newcommand\II{{\mathcal I}}
\newcommand\JJ{{\mathcal J}}
\newcommand\KK{{\mathcal K}}
\newcommand\LL{{\mathcal L}}
\newcommand\MM{{\mathcal M}}
\newcommand\NN{{\mathcal N}}
\newcommand\OO{{\mathcal O}}
\newcommand\PP{{\mathcal P}}
\newcommand\QQ{{\mathcal Q}}
\newcommand\RR{{\mathcal R}}
\newcommand\SSS{{\mathcal S}}
\newcommand\TT{{\mathcal T}}
\newcommand\UU{{\mathcal U}}
\newcommand\VV{{\mathcal V}}
\newcommand\WW{{\mathcal W}}
\newcommand\XX{{\mathcal X}}
\newcommand\YY{{\mathcal Y}}
\newcommand\ZZ{{\mathcal Z}}
\newcommand\CH{{\CC\HH}}
\newcommand\PEY{{\PP\EE\YY}}
\newcommand\MF{{\MM\FF}}
\newcommand\RCT{{{\mathcal R}_{CT}}}
\newcommand\PMF{{\PP\kern-2pt\MM\FF}}
\newcommand\FL{{\FF\LL}}
\newcommand\PML{{\PP\kern-2pt\MM\LL}}
\newcommand\GL{{\GG\LL}}
\newcommand\Pol{{\mathcal P}}
\newcommand\half{{\textstyle{\frac12}}}
\newcommand\Half{{\frac12}}
\newcommand\Mod{\operatorname{Mod}}
\newcommand\Area{\operatorname{Area}}
\newcommand\ep{\epsilon}
\newcommand\hhat{\widehat}
\newcommand\Proj{{\mathbf P}}
\newcommand\U{{\mathbf U}}
 \newcommand\Hyp{{\mathbf H}}
\newcommand\D{{\mathbf D}}
\newcommand\Z{{\mathbb Z}}
\newcommand\R{{\mathbb R}}
\newcommand\Q{{\mathbb Q}}
\newcommand\E{{\mathbb E}}
\newcommand\til{\widetilde}
\newcommand\length{\operatorname{length}}
\newcommand\tr{\operatorname{tr}}
\newcommand\gesim{\succ}
\newcommand\lesim{\prec}
\newcommand\simle{\lesim}
\newcommand\simge{\gesim}
\newcommand{\simmult}{\asymp}
\newcommand{\simadd}{\mathrel{\overset{\text{\tiny $+$}}{\sim}}}
\newcommand{\ssm}{\setminus}
\newcommand{\diam}{\operatorname{diam}}
\newcommand{\pair}[1]{\langle #1\rangle}
\newcommand{\T}{{\mathbf T}}
\newcommand{\inj}{\operatorname{inj}}
\newcommand{\pleat}{\operatorname{\mathbf{pleat}}}
\newcommand{\short}{\operatorname{\mathbf{short}}}
\newcommand{\vertices}{\operatorname{vert}}
\newcommand{\collar}{\operatorname{\mathbf{collar}}}
\newcommand{\bcollar}{\operatorname{\overline{\mathbf{collar}}}}
\newcommand{\I}{{\mathbf I}}
\newcommand{\tprec}{\prec_t}
\newcommand{\fprec}{\prec_f}
\newcommand{\bprec}{\prec_b}
\newcommand{\pprec}{\prec_p}
\newcommand{\ppreceq}{\preceq_p}
\newcommand{\sprec}{\prec_s}
\newcommand{\cpreceq}{\preceq_c}
\newcommand{\cprec}{\prec_c}
\newcommand{\topprec}{\prec_{\rm top}}
\newcommand{\Topprec}{\prec_{\rm TOP}}
\newcommand{\fsub}{\mathrel{\scriptstyle\searrow}}
\newcommand{\bsub}{\mathrel{\scriptstyle\swarrow}}
\newcommand{\fsubd}{\mathrel{{\scriptstyle\searrow}\kern-1ex^d\kern0.5ex}}
\newcommand{\bsubd}{\mathrel{{\scriptstyle\swarrow}\kern-1.6ex^d\kern0.8ex}}
\newcommand{\fsubeq}{\mathrel{\raise-.7ex\hbox{$\overset{\searrow}{=}$}}}
\newcommand{\bsubeq}{\mathrel{\raise-.7ex\hbox{$\overset{\swarrow}{=}$}}}
\newcommand{\tw}{\operatorname{tw}}
\newcommand{\base}{\operatorname{base}}
\newcommand{\trans}{\operatorname{trans}}
\newcommand{\rest}{|_}
\newcommand{\bbar}{\overline}
\newcommand{\UML}{\operatorname{\UU\MM\LL}}
\newcommand{\EL}{\mathcal{EL}}
\newcommand{\tsum}{\sideset{}{'}\sum}
\newcommand{\tsh}[1]{\left\{\kern-.9ex\left\{#1\right\}\kern-.9ex\right\}}
\newcommand{\Tsh}[2]{\tsh{#2}_{#1}}
\newcommand{\qeq}{\mathrel{\approx}}
\newcommand{\Qeq}[1]{\mathrel{\approx_{#1}}}
\newcommand{\qle}{\lesssim}
\newcommand{\Qle}[1]{\mathrel{\lesssim_{#1}}}
\newcommand{\simp}{\operatorname{simp}}
\newcommand{\vsucc}{\operatorname{succ}}
\newcommand{\vpred}{\operatorname{pred}}
\newcommand\fhalf[1]{\overrightarrow {#1}}
\newcommand\bhalf[1]{\overleftarrow {#1}}
\newcommand\sleft{_{\text{left}}}
\newcommand\sright{_{\text{right}}}
\newcommand\sbtop{_{\text{top}}}
\newcommand\sbot{_{\text{bot}}}
\newcommand\sll{_{\mathbf l}}
\newcommand\srr{_{\mathbf r}}
\newcommand\geod{\operatorname{\mathbf g}}
\newcommand\mtorus[1]{\boundary U(#1)}
\newcommand\A{\mathbf A}
\newcommand\Aleft[1]{\A\sleft(#1)}
\newcommand\Aright[1]{\A\sright(#1)}
\newcommand\Atop[1]{\A\sbtop(#1)}
\newcommand\Abot[1]{\A\sbot(#1)}
\newcommand\boundvert{{\boundary_{||}}}
\newcommand\storus[1]{U(#1)}
\newcommand\Momega{\omega_M}
\newcommand\nomega{\omega_\nu}
\newcommand\twist{\operatorname{tw}}
\newcommand\modl{M_\nu}
\newcommand\MT{{\mathbb T}}
\newcommand\Teich{{\mathcal T}}
\renewcommand{\Re}{\operatorname{Re}}
\renewcommand{\Im}{\operatorname{Im}}

\title{Circumcenter extension of Moebius maps to CAT(-1) spaces}

\author{Kingshook Biswas}
\address{Indian Statistical Institute, Kolkata, India. Email: kingshook@isical.ac.in}

\begin{abstract} Given a Moebius homeomorphism $f : \partial X \to \partial Y$
between boundaries of proper, geodesically complete CAT(-1) spaces $X,Y$, we describe an extension
$\hat{f} : X \to Y$ of $f$, called the circumcenter map of $f$, which is constructed using
circumcenters of expanding sets. The extension $\hat{f}$ is shown to coincide with the
$(1, \log 2)$-quasi-isometric extension constructed in \cite{biswas3}, and is locally $1/2$-Holder
continuous. When $X,Y$ are complete, simply connected
manifolds with sectional curvatures $K$ satisfying $-b^2 \leq K \leq -1$ for some $b \geq 1$ then the
extension $\hat{f} : X \to Y$ is a $(1, (1 - \frac{1}{b})\log 2)$-quasi-isometry. Circumcenter extension of Moebius maps is
natural with respect to composition with isometries.
\end{abstract}

\bigskip

\maketitle

\tableofcontents

\section{Introduction}

\medskip

Let $X$ be a CAT(-1) space. There is a positive function called the {\it cross-ratio}
on the space of quadruples of distinct points of the boundary at infinity $\partial X$ of $X$,
defined for $\xi, \xi', \eta, \eta' \in \partial X$ by
$$
[\xi, \xi', \eta, \eta'] = \lim_{a\to\xi, b\to\xi', c\to\eta,d\to\eta'} \exp\left( \frac{1}{2} \left( d(a,c)+d(b,d)-d(a,d)-d(b,c)\right)\right)
$$
(where $a,b,c,d \in X$ converge radially towards $\xi,\xi',\eta,\eta'$).
A map between boundaries of CAT(-1) spaces is called {\it Moebius} if it preserves cross-ratios.
Any isometry between CAT(-1) spaces extends to a Moebius homeomorphism between their boundaries. A
classical fact which turns out to be crucial in many rigidity results including the Mostow Rigidity
theorem is that a Moebius map from the boundary of real hyperbolic space to itself extends to an isometry.
More generally Bourdon showed (\cite{bourdon2}) that if $X$ is a rank one symmetric
space of noncompact type with maximum of sectional curvatures equal to -1 and $Y$ a CAT(-1) space then any Moebius embedding $f :
\partial X \to \partial Y$ extends to an isometric embedding $F : X \to Y$. In \cite{biswas3} the problem
of extending Moebius maps was considered for general CAT(-1) spaces, where it was shown that any Moebius homeomorphism
$f : \partial X \to \partial Y$ between boundaries of proper, geodesically complete CAT(-1) spaces $X, Y$ extends
to a $(1, \log 2)$-quasi-isometry $F : X \to Y$. The proof of this theorem uses an
isometric embedding of a proper, geodesically complete CAT(-1) space into a certain space of Moebius metrics on the boundary of
the space. A nearest point projection to the subspace of visual metrics is used to construct the extension.
We show that this nearest point is unique, and can be constructed as a limit of circumcenters of certain expanding sets. The
extension constructed in \cite{biswas3} is thus uniquely determined.
We call the extension the circumcenter map of $f$. It is readily seen to satisfy naturality properties with
respect to composition with isometries. We have:

\medskip

\begin{theorem} \label{thmone} Let $X, Y$ be proper, geodesically complete CAT(-1) spaces, and $f : \partial X \to \partial Y$
a Moebius homeomorphism. Then the circumcenter extension $\hat{f} : X \to Y$ of $f$ is a $(1, \log 2)$-quasi-isometry
which is locally $1/2$-Holder continuous: 
$$
d(\hat{f}(x), \hat{f}(y)) \leq 2 d(x,y)^{1/2}
$$
for all $x,y \in X$ such that $d(x,y) \leq 1$.
\end{theorem}

\medskip

When the spaces $X, Y$ are also assumed to be manifolds with curvature bounded below
 we have the following improvement on the main result of \cite{biswas3}:

\medskip

\begin{theorem} \label{mainthm} Let $X, Y$ be complete, simply connected Riemannian manifolds with sectional curvatures
satisfying $-b^2 \leq K \leq -1$ for some constant $b \geq 1$. For any Moebius homeomorphism $f : \partial X \to \partial Y$, the
circumcenter extension $\hat{f}$ of $f$ is a $(1, (1 - \frac{1}{b})\log 2)$-quasi-isometry $\hat{f} : X \to Y$ with image
$\frac{1}{2}(1 - \frac{1}{b})\log 2$-dense
in $Y$.
\end{theorem}

\medskip

\medskip

We mention that one of the motivations for considering the problem of extending Moebius
maps is the {\it marked length spectrum rigidity} problem. This asks whether an isomorphism
$\phi : \pi_1(X) \to \pi_1(Y)$ between fundamental groups of closed negatively curved
manifolds which preserves lengths of closed geodesics (recall that in negative curvature
each homotopy class of closed curves contains a unique closed geodesic) is necessarily induced
by an isometry $F : X \to Y$. Otal (\cite{otal2}) proved that this is indeed the case in dimension two.
The problem remains open in higher dimensions. It is known however
to be equivalent to the {\it geodesic conjugacy} problem, which asks whether the
existence of a homeomorphism between the unit tangent bundles
$\phi : T^1 X \to T^1 Y$ conjugating the geodesic flows
implies isometry of the manifolds. Hamenstadt (\cite{hamenstadt1}) proved that equality
of marked length spectra is equivalent to existence of a
geodesic conjugacy.

\medskip

Bourdon showed in \cite{bourdon1}, that for a
Gromov-hyperbolic group $\Gamma$ with two quasi-convex actions on
CAT(-1) spaces $X, Y$, the natural $\Gamma$-equivariant
homeomorphism $f$ between the limit sets $\Lambda X, \Lambda Y$ is Moebius if and only if
there is a $\Gamma$-equivariant conjugacy of the abstract geodesic
flows $\mathcal{G}\Lambda X$ and $\mathcal{G} \Lambda Y$ compatible with
$f$. In particular for $\tilde{X}, \tilde{Y}$ the universal covers of two closed
negatively curved manifolds $X, Y$ (with sectional curvatures bounded above by $-1$), the geodesic flows of $X, Y$ are topologically conjugate if
and only if the induced equivariant boundary map $f : \partial X
\to \partial Y$ is Moebius. Thus an affirmative answer to the problem of extending Moebius maps to
isometries would also yield a solution to the equivalent problems of marked length spectrum rigidity and
geodesic conjugacy.

\medskip

Finally we remark that in \cite{biswas4} it is proved that in certain cases Moebius maps
between boundaries of simply connected negatively curved manifolds do extend to isometries
(more precisely, local and infinitesimal rigidity results are proved for deformations of the metric
on a compact set).

\medskip

\section{Spaces of Moebius metrics}

\medskip

We recall in this section the definitions and facts from \cite{biswas3} which we will be needing.

\medskip

Let $(Z,\rho_0)$ be a compact metric space with at least four points. For a metric $\rho$ on
$Z$ we define the metric cross-ratio with respect to $\rho$ of a quadruple of distinct
points $(\xi, \xi', \eta, \eta')$ of $Z$ by
$$
[\xi \xi' \eta \eta']_{\rho} := \frac{\rho(\xi, \eta) \rho(\xi', \eta')}{\rho(\xi,
\eta')\rho(\xi', \eta)}
$$
We say that a diameter one metric $\rho$ on $Z$ is {\it antipodal} if for
any $\xi \in Z$ there exists $\eta \in Z$ such that $\rho(\xi,
\eta) = 1$. We assume that $\rho_0$ is diameter one and antipodal. We say two metrics $\rho_1, \rho_2$
on $Z$ are {\it Moebius equivalent} if their metric cross-ratios agree:
$$
[\xi \xi' \eta \eta']_{\rho_1} = [\xi \xi' \eta \eta']_{\rho_2}
$$
for all $(\xi, \xi', \eta, \eta')$. The space of Moebius metrics on $Z$ is defined to be
$$
\mathcal{M}(Z, \rho_0) := \{ \rho : \rho \hbox{ is an antipodal, diameter one metric on } Z
\hbox{ Moebius equivalent to } \rho_0 \}
$$
We will write $\mathcal{M}(Z, \rho_0) = \mathcal{M}$. We have the following from \cite{biswas3}:

%
\medskip

\begin{theorem} \label{derivative} For any $\rho_1, \rho_2 \in \mathcal{M}$, there is a positive continuous function
$\frac{d\rho_2}{d\rho_1}$ on $Z$, called the derivative of $\rho_2$ with respect to $\rho_1$, such that
the following holds (the "Geometric Mean Value Theorem"):
$$
\rho_2(\xi, \eta)^2 = \frac{d\rho_2}{d\rho_1}(\xi) \frac{d\rho_2}{d\rho_1}(\eta) \rho_1(\xi, \eta)^2
$$
for all $\xi, \eta \in Z$.

\smallskip

Moreover for $\rho_1,
\rho_2, \rho_3 \in \mathcal{M}$ we have
$$
\frac{d\rho_3}{d\rho_1} = \frac{d\rho_3}{d\rho_2} \frac{d\rho_2}{d\rho_1}
$$
and
$$
\frac{d\rho_2}{d\rho_1} = 1/\left(\frac{d\rho_1}{d\rho_2}\right)
$$
\end{theorem}

\medskip

\begin{lemma} \label{maxmin}
$$
\max_{\xi \in Z} \frac{d\rho_2}{d\rho_1}(\xi) \cdot \min_{\xi \in Z}
\frac{d\rho_2}{d\rho_1}(\xi) = 1
$$
Moreover if $\frac{d\rho_2}{d\rho_1}$ attains its maximum 
at $\xi$ and $\rho_1(\xi, \eta) = 1$ then $\frac{d\rho_2}{d\rho_1}$ attains its minimum
at $\eta$, and $\rho_2(\xi, \eta) = 1$.
\end{lemma}

\medskip

\noindent{\bf Proof:} Let $\lambda, \mu$ denote the maximum and
minimum values of $\frac{d\rho_2}{d\rho_1}$ respectively, and let $\xi, \xi' \in Z$ denote
points where the maximum and minimum values are attained
respectively. Given $\eta \in Z$ such that $\rho_1(\xi, \eta) =
1$, we have, using the Geometric Mean-Value Theorem,
$$
1 \geq \rho_2(\xi, \eta)^2 = \frac{d\rho_2}{d\rho_1}(\xi)
\frac{d\rho_2}{d\rho_1}(\eta) \geq \lambda \cdot
\mu
$$
while choosing $\eta' \in Z$ such that $\rho_2(\xi', \eta') = 1$, we have
$$
1 \geq \rho_1(\xi', \eta')^2 =
1/\left(\frac{d\rho_2}{d\rho_1}(\xi')
\frac{d\rho_2}{d\rho_1}(\eta')\right) \geq 1/(\lambda
\mu)
$$
hence $\lambda \cdot \mu = 1$.

\medskip

By the above we have
$$
\frac{d\rho_2}{d\rho_1}(\eta) \leq 1/\frac{d\rho_2}{d\rho_1}(\xi) = 1/\lambda = \mu
$$
hence $\frac{d\rho_2}{d\rho_1}(\eta) = \mu$. By the Geometric Mean Value Theorem this gives
$$
\rho_2(\xi ,\eta)^2 = \rho_1(\xi, \eta)^2 \frac{d\rho_2}{d\rho_1}(\xi) \frac{d\rho_2}{d\rho_1}(\eta) = 1 \cdot \lambda \cdot \mu = 1
$$
$\diamond$

\medskip

For $\rho_1, \rho_2 \in \mathcal{M}$, we define
$$
d_{\mathcal{M}}(\rho_1, \rho_2) := \max_{\xi \in Z}
\log \frac{d\rho_2}{d\rho_1}(\xi)
$$

\medskip

From \cite{biswas3} we have:

\medskip

\begin{lemma} \label{metric} The function $d_{\mathcal{M}}$
defines a metric on $\mathcal{M}$. The metric space $(\mathcal{M}, d_{\mathcal{M}})$ is proper.
\end{lemma}

\medskip

\bigskip

\section{Visual metrics on the boundary of a CAT(-1) space}

\medskip

Let $X$ be a proper CAT(-1) space such that $\partial X$ has at least four points.

\medskip

We recall below the definitions and some
elementary properties of visual metrics and Busemann functions;
for proofs we refer to \cite{bourdon1}:

\medskip

Let $x \in X$ be a basepoint. The {\it Gromov product} of two
points $\xi, \xi' \in \partial X$ with respect to $x$ is defined by
$$
(\xi | \xi')_x = \lim_{(a,a') \to (\xi, \xi')}
\frac{1}{2}(d(x,a) + d(x,a') - d(a,a'))
$$
where $a,a'$ are points of $X$ which converge radially towards
$\xi$ and $\xi'$ respectively. The {\it visual metric} on
$\partial X$ based at the point $x$ is defined by
$$
\rho_x(\xi, \xi') := e^{-(\xi|\xi')_x}
$$
The distance $\rho_x(\xi,\xi')$ is less than or equal to one, with
equality iff $x$ belongs to the geodesic $(\xi \xi')$.

\medskip

\begin{lemma} \label{visualantipodal} If $X$ is geodesically
complete then $\rho_x$ is a diameter one antipodal metric.
\end{lemma}

\medskip

The Busemann function $B : \partial X \times X \times X \to
\mathbb{R}$ is defined by
$$
B(x, y, \xi) := \lim_{a \to \xi} d(x,a) - d(y,a)
$$
where $a \in X$ converges radially towards $\xi$.

\medskip

\begin{lemma} \label{busemann} We have $|B(x,y,\xi)| \leq
d(x,y)$ for all $\xi \in \partial X, x,y \in X$. Moreover
$B(x,y,\xi) = d(x,y)$ iff $y$ lies on the geodesic ray $[x,
\xi)$ while $B(x,y,\xi) = -d(x,y)$ iff $x$ lies on the
geodesic ray $[y, \xi)$.
\end{lemma}

\medskip

We recall the following Lemma from \cite{bourdon1}:

\medskip

\begin{lemma} \label{visualmvt} For $x, y \in X, \xi, \eta \in
\partial X$ we have
$$
\rho_y(\xi, \eta)^2 = \rho_x(\xi, \eta)^2 e^{B(x,y,\xi)} e^{B(x,y,\eta)}
$$
\end{lemma}

\medskip

An immediate corollary of the above Lemma is the following:

\medskip

\begin{lemma} \label{visualmoebius} The visual metrics $\rho_x, x \in X$
are Moebius equivalent to each other and
$$
\frac{d\rho_y}{d\rho_x}(\xi) = e^{B(x,y,\xi)}
$$
\end{lemma}

\medskip

It follows that the metric cross-ratio $[\xi\xi'\eta\eta']_{\rho_x}$ of a quadruple
$(\xi, \xi',\eta,\eta')$ is independent of the choice of $x \in
X$. Denoting this common value by $[\xi\xi'\eta\eta']$, it is
shown in \cite{bourdon2} that the cross-ratio is given by
$$
[\xi\xi'\eta\eta'] = \lim_{(a,a',b,b') \to (\xi, \xi',\eta,\eta')} \exp(\frac{1}{2}(d(a,b)+d(a',b') -
d(a,b') - d(a',b)))
$$
where the points $a,a',b,b' \in X$ converge radially towards
$\xi,\xi',\eta,\eta' \in \partial X$.

\medskip

We assume henceforth that $X$ is a proper, geodesically complete CAT(-1) space.
We let $\mathcal{M} = \mathcal{M}(\partial X, \rho_x)$ (this space is independent of the
choice of $x \in X$). From \cite{biswas3} we have:

\medskip

\begin{lemma} The map
\begin{align*}
i_X : X & \to \mathcal{M} \\
         x & \mapsto \rho_x  \\
\end{align*}
is an isometric embedding and the image is closed in $\mathcal{M}$.
\end{lemma}

\medskip

For $k > 0$ and $y,z \in X$ distinct from $x \in X$ let $\angle^{(-k^2)} y x z \in [0, \pi]$ denote the
angle at the vertex $\overline{x}$ in a comparison triangle $\overline{x}\overline{y}\overline{z}$ in the model space $\mathbb{H}_{-k^2}$ of
constant curvature $-k^2$.

\medskip

\begin{lemma} \label{compexist} For $\xi, \eta \in \partial X$, the limit of the comparison angles $\angle^{(-k^2)} y x z$ exists as $y,z$ converge to $\xi, \eta$ along
the geodesic rays $[x,\xi), [x, \eta)$ respectively. Denoting this limit by $\angle^{(-k^2)} \xi x \eta$, it satisfies
$$
\sin \left( \frac{\angle^{(-k^2)} \xi x \eta}{2} \right) = \rho_x(\xi, \eta)^k
$$
\end{lemma}

\medskip

\noindent{\bf Proof:} A comparison triangle in $\mathbb{H}_{-k^2}$ with side lengths $a = d(x,y), b = d(x, z), c = d(y,z)$ and angle
$\theta = \angle^{(-k^2)} y x z$ at the vertex corresponding to $x$ corresponds to a triangle in $\mathbb{H}_{-1}$ with side lengths
$ka, kb, kc$ and angle $\theta$ at the vertex opposite the side with length $kc$. By the hyperbolic law of cosine we have
$$
\cosh kc = \cosh ka \cosh kb - \sinh ka \sinh kb \cos \theta
$$
As $y \to \xi, z \to \eta$, we have $a, b, c \to \infty$, and $a + b - c \to 2(\xi|\eta)_x$, thus
\begin{align*}
\cos \theta & = \frac{\cosh ka \cosh kb}{\sinh ka \sinh kb} - \frac{\cosh kc}{\sinh ka \sinh kb} \\
            & \to 1 - 2 e^{-2k(\xi|\eta)_x} \\
\end{align*}
hence the angle $\theta$ converges to a limit. Denoting this limit by $\angle^{(-k^2)} \xi x \eta$, by the above
it satisfies
$$
\cos(\angle^{(-k^2)} \xi x \eta) = 1 - 2\rho_x(\xi, \eta)^{2k}
$$
and hence
$$
\sin \left( \frac{\angle^{(-k^2)} \xi x \eta}{2} \right) = \rho_x(\xi, \eta)^k
$$
$\diamond$

\medskip

\begin{lemma} \label{buseform} For $x,y \in X, \xi \in \partial X$ and $k > 0$,
the limit of the comparison angles $\angle^{(-k^2)} y x z$ exists as $z$ converges to $\xi$ along
the geodesic ray $[x,\xi)$. Denoting this limit by $\angle^{(-k^2)} y x \xi$, it satisfies
$$
e^{kB(y,x,\xi)} = \cosh(kd(x,y)) - \sinh(kd(x,y)) \cos(\angle^{(-k^2)}y x \xi)
$$
\end{lemma}

\medskip

\noindent{\bf Proof:} A comparison triangle in $\mathbb{H}_{-k^2}$ with side lengths $a = d(x,y), b = d(x, z), c = d(y,z)$ and angle
$\theta = \angle^{(-k^2)} y x z$ at the vertex corresponding to $x$ corresponds to a triangle in $\mathbb{H}_{-1}$ with side lengths
$ka, kb, kc$ and angle $\theta$ at the vertex opposite the side with length $kc$. By the hyperbolic law of cosine we have
$$
\cosh kc = \cosh ka \cosh kb - \sinh ka \sinh kb \cos \theta
$$
As $z \to \xi$, we have $b, c \to \infty$, and $c - b \to B(y, x, \xi)$, thus
\begin{align*}
\cos \theta & = \frac{\cosh ka \cosh kb}{\sinh ka \sinh kb} - \frac{\cosh kc}{\sinh ka \sinh kb} \\
            & \to \frac{\cosh ka}{\sinh ka} - \frac{e^{kB(y,x,\xi)}}{\sinh ka} \\
\end{align*}
hence the angle $\theta$ converges to a limit. Denoting this limit by $\angle^{(-k^2)} \xi x \eta$, by the above
it satisfies
$$
e^{kB(y,x,\xi)} = \cosh(kd(x,y)) - \sinh(kd(x,y)) \cos(\angle^{(-k^2)}y x \xi)
$$
$\diamond$

\medskip

We now consider the behaviour of the derivatives $\frac{d\rho_y}{d\rho_x}$ as $t = d(x,y) \to 0$ and the point $y$
converges radially towards $x$ along a geodesic. For functions $F_t$ on $\partial X$ we write $F_t = o(t)$ if
$||F_t||_{\infty} = o(t)$. We have the following formula from \cite{biswas3}, which may be thought of as a formula
for the derivative of the map $i_X$ along a geodesic:

\medskip

\begin{lemma} \label{embedderiv} As $t \to 0$ we have
$$
\log \frac{d\rho_y}{d\rho_x}(\xi) = t \cos(\angle^{(-1)}y x \xi) + o(t)
$$
\end{lemma}

\medskip

\section{Conformal maps, Moebius maps and geodesic conjugacies}

\medskip

We start by recalling the definitions of conformal maps, Moebius
maps, and the abstract geodesic flow of a CAT(-1) space.

\medskip

\begin{definition} A homeomorphism between metric spaces $f :
(Z_1, \rho_1) \to (Z_2, \rho_2)$ with no isolated points is said to be {\it conformal} if
for all $\xi \in Z_1$, the limit
$$
df_{\rho_1, \rho_2}(\xi) := \lim_{\eta \to \xi} \frac{\rho_2(f(\xi),
f(\eta))}{\rho_1(\xi, \eta)}
$$
exists and is positive. The positive function $df_{\rho_1,
\rho_2}$ is called the derivative of $f$ with respect to $\rho_1, \rho_2$.
We say $f$ is {\it $C^1$ conformal} if its derivative is continuous.

\medskip

Two metrics $\rho_1, \rho_2$ inducing the same topology on a set
$Z$, such that $Z$ has no isolated points,
are said to be conformal (respectively $C^1$ conformal) if the
map $id_Z : (Z, \rho_1) \to (Z, \rho_2)$ is conformal
(respectively $C^1$ conformal). In this case we denote the
derivative of the identity map by $\frac{d\rho_2}{d\rho_1}$.
\end{definition}

\medskip

\begin{definition} A homeomorphism between metric spaces $f :
(Z_1, \rho_1) \to (Z_2, \rho_2)$ (where $Z_1$ has at least four
points) is said to be Moebius if it preserves metric cross-ratios
with respect to $\rho_1, \rho_2$. The derivative of $f$ is defined
to be the derivative $\frac{df_*\rho_2}{\rho_1}$ of the Moebius
equivalent metrics $f_* \rho_2, \rho_1$ as defined in section 2
(where $f_* \rho_2$ is the pull-back of $\rho_2$ under $f$).
\end{definition}

\medskip

From the results of section 2 it follows that any Moebius map
between compact metric spaces with no isolated points is $C^1$ conformal,
and the two definitions of the derivative of $f$ given above
coincide. Moreover any Moebius map $f$ satisfies the geometric
mean-value theorem,
$$
\rho_2(f(\xi), f(\eta))^2 = \rho_1(\xi,\eta)^2
df_{\rho_1,\rho_2}(\xi) df_{\rho_1,\rho_2}(\xi)
$$

\medskip

\begin{definition} Let $(X, d)$ be a CAT(-1) space. The abstract geodesic flow
space of $X$ is defined to be the space of bi-infinite geodesics
in $X$,
$$
\mathcal{G}X := \{ \gamma : (-\infty,+\infty) \to X | \gamma
\hbox{ is an isometric embedding} \}
$$
endowed with the topology of uniform convergence on compact
subsets. This topology is metrizable with a distance defined by
$$
d_{\mathcal{G}X}(\gamma_1, \gamma_2):= \int_{-\infty}^{\infty}
d(\gamma_1(t), \gamma_2(t)) \frac{e^{-|t|}}{2} \ dt
$$
We define also two projections
\begin{align*}
\pi : \mathcal{G}X & \to X \\
             \gamma  & \mapsto \gamma(0) \\
\end{align*}
and
\begin{align*}
p : \mathcal{G}X & \to \partial X \\
             \gamma  & \mapsto \gamma(+\infty) \\
\end{align*}

It is shown in Bourdon \cite{bourdon1} that $\pi$ is
$1$-Lipschitz.

\medskip

For $x \in X$, the unit tangent sphere $T^1_x X \subset \mathcal{G}X$ is defined to be
$$
T^1_x X := \pi^{-1}(x)
$$

\medskip

The abstract geodesic flow of $X$ is defined to be the one-parameter group
of homeomorphisms
\begin{align*}
\phi_t : \mathcal{G}X & \to \mathcal{G}X \\
         \gamma       & \mapsto \gamma_t \\
\end{align*}
for $t \in \mathbb{R}$, where $\gamma_t$ is the geodesic $s
\mapsto \gamma(s+t)$.

\medskip

The flip is defined to be the map
\begin{align*}
\mathcal{F} : \mathcal{G}X & \to \mathcal{G}X \\
              \gamma & \mapsto \overline{\gamma} \\
\end{align*}
where $\overline{\gamma}$ is the geodesic $s
\mapsto \gamma(-s)$.
\end{definition}

\medskip

We observe that for a simply connected complete Riemannian manifold $X$ with
sectional curvatures bounded above by $-1$, the map
\begin{align*}
\mathcal{G}X & \to T^1 X \\
 \gamma      & \mapsto \gamma'(0) \\
\end{align*}
is a homeomorphism conjugating the abstract geodesic flow of $X$
to the usual geodesic flow of $X$ and the flip $\mathcal{F}$ to the
usual flip on $T^1 X$.

\medskip

Let $f : \partial X \to \partial Y$ be a
conformal map between the boundaries of CAT(-1) spaces $X, Y$
equipped with visual metrics. Then $f$ induces a bijection $\phi_f :
\mathcal{G}X \to \mathcal{G}Y$ conjugating the geodesic flows, which is defined as
follows:

\medskip

Given $\gamma \in \mathcal{G}X$, let
$\gamma(-\infty) = \xi, \gamma(+\infty) = \eta, x = \gamma(0)$, then there
is a unique point $y$ on the bi-infinite geodesic $(f(\xi),f(\eta))$ such that $df_{\rho_x,
\rho_y}(\eta) = 1$. Define $\phi_f(\gamma) = \gamma^*$ where
$\gamma^*$ is the unique geodesic in $Y$ satisfying $\gamma^*(-\infty) = f(\xi),
\gamma^*(+\infty) = f(\eta), \gamma^*(0) = y$. Then $\phi_f : \mathcal{G}X \to
\mathcal{G}Y$ is a bijection conjugating the geodesic flows. From \cite{biswas3} we
have:

\medskip

\begin{prop} \label{confconj} The map $\phi_f$ is a homeomorphism if $f$ is $C^1$ conformal. If $f$ is
Moebius then $\phi_f$ is flip-equivariant.
\end{prop}

\medskip

\section{Circumcenters of expanding sets and $\mathcal{F}K$-convex functions}

\medskip

Let $X$ be a proper, geodesically complete CAT(-1) space. Recall that for any bounded subset $B$ of
 $X$, there is a unique point $x$
which minimizes the function
$$
z \mapsto \sup_{y \in B} d(z, y)
$$
The point $x$ is called the {\it circumcenter} of $B$, and the number $\sup_{y \in B} d(x,y)$ is called
the circumradius of B. We will denote these by $c(B)$ and $r(B)$ respectively.

\medskip

Given $K \leq 0$, a function $f : X \to \bbR$ is said to be $\mathcal{F}K$-convex if it is continuous and its restriction to any
geodesic satisfies $f'' + Kf \geq 0$ in the barrier sense. This means that $f \leq g$ if $g$ coincides with
$f$ at the endpoints of a subsegment and satisfies $g'' + Kg = 0$.
We have the following from \cite{alexanderbishop}:

\medskip

\begin{prop} Let $y \in X, \xi \in \partial X$. Then:

\smallskip

\noindent (1) The function $x \mapsto \cosh(d(x,y))$ is $\mathcal{F}(-1)$-convex.

\smallskip

\noindent (2) The function $x \mapsto \exp(B(x,y,\xi))$ is $\mathcal{F}(-1)$-convex.

\end{prop}

\medskip

\begin{prop} \label{fkmin} Let $f$ be a positive, proper, $\mathcal{F}(-1)$-convex function on $X$.
Then $f$ attains its minimum at a unique point $x \in X$.
\end{prop}

\medskip

\noindent{\bf Proof:} Since $f$ is continuous, bounded below, and proper, $f$ attains its minimum
at some $x \in X$. If $x' \neq x$ is another point where $f$ attains its minimum, let
$\gamma : [-d, d] \to X$ be the geodesic joining $x$ to $x'$ where $d = d(x,x')/2 > 0$.
Then $g(t) = f(x) \cosh t / \cosh d$ satisfies
$g'' - g = 0$, and agrees with $f$ at the endpoints of $\gamma$, hence $f(\gamma(0)) \leq g(0) = f(x) / \cosh d < f(x)$,
a contradiction. $\diamond$

\medskip

\begin{prop} \label{minconv} Let $f_n, f$ be positive, proper, $\mathcal{F}(-1)$-convex functions on $X$
such that $f_n \to f$ uniformly on compacts. If $x_n, x$ denote the points where $f_n, f$ attain their
minima, then $x_n \to x$.
\end{prop}

\medskip

\noindent{\bf Proof:} We first show that $\{x_n\}$ is bounded. If not, passing to a subsequence
we may assume $d(x, x_n) \to +\infty$. For $n$ sufficiently large we have $f_n(x) \leq 2f(x)$. Thus
$f_n(x_n) \leq f_n(x) \leq 2f(x)$ as well. Let $\gamma_n : [-d_n, d_n] \to X$ be the unique
geodesic joining $x$ to $x_n$, where $d_n = d(x, x_n)/2$. Then the function
$$
g(t) = \frac{1}{\sinh(2d_n)}\left[ (\sinh d_n)(f_n(x)+f_n(x_n)) \cosh t + (\cosh d_n)(f_n(x_n) - f_n(x)) \sinh t \right]
$$
satisfies $g'' - g = 0$, and agrees with $f_n$ the endpoints of $\gamma_n$. Thus for any $s > 0$, for $n$ large such that
$s < d_n$, letting $y_n = \gamma_n(-d_n + s)$, we have
\begin{align*}
f_n(y_n) & \leq g(-d_n + s) \\
         &\leq \frac{1}{2 \sinh d_n \cosh d_n} \left[ (\sinh d_n \cosh (d_n - s))(4f(x)) + \cosh d_n \sinh (d_n - s) (4f(x)) \right] \\
\end{align*}

Since $d_n \to +\infty$ this implies that for $n$ sufficiently large we have
$$
f_n(y_n) \leq \frac{1}{2} e^{-s}(1 + o(1)) (8f(x)) < f(x)/2
$$
for $s > 0$ large enough. Fixing such an $s$, the points $y_n = \gamma_n(-d_n+s)$ lie in the closed ball
$B$ of radius $s$ around $x$, so passing to a subsequence we may assume that $y_n \to y \in B$. Since $f_n \to f$
uniformly on $B$, $f_n(y_n) \to f(y)$, hence $f(y) \leq f(x)/2 < f(x)$, a contradiction.

\medskip

Thus the sequence $\{ x_n \}$ is bounded. To show $x_n \to x$, it suffices to show that the only limit point
of $\{x_n\}$ is $x$. Let $K$ be a compact containing $\{x_n\}$. Suppose $x_{n_k} \to y$. Then $f_{n_k}(x_{n_k}) \leq f_{n_k}(x)$ for all
$k$. Since $f_{n_k} \to f$ uniformly on $K$, letting $k$ tend to infinity gives $f(y) \leq f(x)$. By the previous proposition this implies
$y = x$. $\diamond$

\medskip

Let $K$ be a compact subset of $\mathcal{G}X$. Define the function
$$
u_K(z) = \sup_{\gamma \in K} \exp(B(z, \pi(\gamma), \gamma(+\infty)))
$$

\begin{prop} \label{fkconvex} The function $u_K$ is a positive, $\mathcal{F}(-1)$-convex function.
It is proper if $p(K) \subset \partial X$ is not a singleton.
\end{prop}

\medskip

\noindent{\bf Proof:} For each $\gamma \in K$, the function $z \mapsto \exp(B(z, \pi(\gamma), \gamma(+\infty)))$
is $\mathcal{F}(-1)$-convex. Thus $u_K$, being the supremum of a family of $\mathcal{F}(-1)$-convex functions,
satisfies the $\mathcal{F}(-1)$-convexity inequality. It remains to show that $u_K$ is continuous.

\medskip

Let $z_n \to z$ in $X$. Define functions $h_n, h : K \to \bbR$ by
$$
h_n(\gamma) := B(z_n, \pi(\gamma), \gamma(+\infty)) , h(\gamma) := B(z, \pi(\gamma), \gamma(+\infty))
$$
Then $|h_n(\gamma) - h(\gamma)| = |B(z_n, z, \gamma(\infty)| \leq d(z_n, z)$, so $h_n \to h$ uniformly on $K$. It
follows that
$$
u_K(z_n) = ||e^{h_n}||_{\infty} \to ||e^{h}||_{\infty} = u_K(z)
$$
Thus $u_K$ is continuous.

\medskip

Now suppose $p(K)$ is not a singleton, so 
there exist $\gamma_1, \gamma_2 \in K$ such that the endpoints $\xi_i = \gamma_i(+\infty), i = 1,2$ are distinct.
Let $x_n$ be a sequence in $X$ tending to infinity. Suppose $u_K(x_n)$ does not tend to $+\infty$. Passing to a
subsequence we may assume $u_K(x_n) \leq M$ for all $n$ for some $M > 0$.
Passing to a further subsequence we may assume $x_n \to \xi \in \partial X$. We can choose a $\xi_i \neq \xi$. 
Let $x = \pi(\gamma_i)$, then by Lemma \ref{buseform} we have
\begin{align*}
\exp(B(x_n, x, \xi_i)) & = \cosh(d(x_n, x)) - \sinh(d(x_n, x)) \cos(\angle^{(-1)}x_n x \xi_i) \\
                       & = e^{-d(x_n, x)} + 2 \sinh(d(x_n, x)) \sin^2\left(\frac{\angle^{(-1)}x_n x \xi_i}{2}\right) \\
                       & \to +\infty \\
\end{align*}
since $\angle^{(-1)}x_n x \xi_i \to \angle^{(-1)}\xi x \xi_i > 0$.   
Hence $u_K(x_n) \geq \exp(B(x_n, x, \xi_i)) \to +\infty$,
a contradiction. This shows that $u_K$ is proper. $\diamond$

\medskip

\begin{definition} Let $K$ be a compact subset of $\mathcal{G}X$ such that $p(K) \subset \partial X$ is not a 
singleton. The asymptotic circumcenter of $K$ is defined to be the unique $x$ in $X$
where the function $u_K$ attains its minimum. We denote the asymptotic circumcenter by $x = c_{\infty}(K)$.
\end{definition}

\medskip

The reason for the name 'asymptotic circumcenter' is explained by the following proposition:

\medskip

\begin{prop} Let $K$ be a compact subset of $\mathcal{G}X$ such that $p(K)$ is not a singleton. 
Define for $t > 0$ bounded subsets $A_t$ of $X$ by $A_t = \pi(\phi_t(K))$,
where $\phi_t$ denotes the geodesic flow on $\mathcal{G}X$. Then
$$
c(A_t) \to c_{\infty}(K)
$$
as $t \to +\infty$, i.e. the circumcenters of the sets $A_t$ converge to the asymptotic circumcenter of $K$.
\end{prop}

\medskip

\noindent{\bf Proof:} Let $u = u_K$, and for $t > 0$ define $u_t : X \to \bbR$ by
$$
u_t(z) = (\sup_{y \in A_t} \cosh(d(z,y))) \cdot 2e^{-t}
$$
It is easy to see that $u_t$ is a positive, proper, $\mathcal{F}(-1)$-convex function, and that the circumcenter
of $A_t$ is the unique minimizer of the function $u_t$. Since $c_{\infty}(K)$ is the unique minimizer of $u$, by the
previous proposition it suffices to show that $u_t \to u$ uniformly on compacts as $t \to \infty$.

\medskip

Note
$$
u_t(z) = (\sup_{\gamma \in K} \cosh(d(z, \gamma(t)))) \cdot 2 e^{-t}
$$
Now for $z$ in a compact ball $B$ and $\gamma$ in the compact $K$,
$$
d(z, \gamma(t)) - t \to B(z, \pi(\gamma), \gamma(+\infty))
$$
as $t \to +\infty$ uniformly in $z \in B, \gamma \in K$. It follows that
$$
\cosh(d(z, \gamma(t))) \cdot 2 e^{-t} \to \exp(B(z, \pi(\gamma), \gamma(+\infty)))
$$
as $t \to +\infty$ uniformly in $z \in B, \gamma \in K$. Since the convergence in $z, \gamma$ is uniform,
the supremums over $\gamma \in K$ converge, uniformly for $z \in B$:
$$
u_t(z) = (\sup_{\gamma \in K} \cosh(d(z, \gamma(t)))) \cdot 2 e^{-t} \to \sup_{\gamma \in K} \exp(B(z, \pi(\gamma), \gamma(+\infty))) = u(z)
$$
uniformly in $z \in B$.

$\diamond$

\bigskip

\section{Circumcenter extension of Moebius maps and nearest point projections}

\medskip

Let $f : \partial X \to \partial Y$ be a Moebius homeomorphism between boundaries of
proper, geodesically complete CAT(-1) spaces $X, Y$, and let $\phi_f : \mathcal{G}X \to \mathcal{G}Y$
denote the associated geodesic conjugacy.

\medskip

\begin{definition} The circumcenter extension of the Moebius map $f$ is the map $\hat{f} : X \to Y$
defined by
$$
\hat{f}(x) := c_{\infty}(\phi_f(T^1_x X)) \in Y
$$
(note that $p(\phi_f(T^1_x X)) = \partial Y$ is not a singleton, so the asymptotic circumcenter of 
$\phi_f(T^1_x X) \subset \mathcal{G}Y$ exists).
\end{definition}

\medskip

In \cite{biswas3}, a $(1, \log 2)$-quasi-isometric extension $F : X \to Y$ of the Moebius map $f$ is constructed as follows.
Since $f$ is Moebius, push-forward by $f$ of metrics on $\partial X$ to metrics on $\partial Y$ gives a
map between the spaces of Moebius metrics $f_* : \mathcal{M}(\partial X) \to \mathcal{M}(\partial Y)$, which is
easily seen to be an isometry. For each $\rho \in \mathcal{M}(\partial Y)$, we can choose a nearest point to $\rho$
in the subspace of visual metrics $i_Y(Y) \subset \mathcal{M}(\partial Y)$. This defines a nearest-point projection
$r_Y : \mathcal{M}(\partial Y) \to Y$. The extension $F$ is then defined by
$$
F = r_Y \circ f_* \circ i_X
$$

\medskip

We show below that if $\rho \in \mathcal{M}(\partial Y)$ is the push-forward of a visual metric on $\partial X$,
$\rho = f_* \rho_x$ for some $x \in X$, then in fact there is a unique visual metric $\rho_y \in \mathcal{M}(\partial Y)$
nearest to $\rho$, given by $y = \hat{f}(x)$, the asymptotic circumcenter of $\phi_f(T^1_x X)$. It follows that the
extension $F$ defined above is uniquely determined and equals the circumcenter extension $\hat{f}$.

\medskip

\begin{prop} \label{nearuniq} Let $x \in X$ and let $\rho = f_* \rho_x \in \mathcal{M}(\partial Y)$. Then $y = \hat{f}(x)$ is the unique
minimizer of the function $z \in Y \mapsto d_{\mathcal{M}}(\rho, \rho_z)$. In particular, $\hat{f} = F$, so $\hat{f}$ is a
$(1, \log 2)$-quasi-isometry.
\end{prop}

\medskip

\noindent{\bf Proof:} Fix a $z \in Y$.
Given $\xi \in \partial X$, let $\gamma \in T^1_x X$ be such that $\gamma(+\infty) = \xi$. Let
$p = \pi(\phi_f(\gamma)) \in Y$. Then by definition of $\phi_f$, we have
$$
\frac{d f_* \rho_x}{d\rho_p}(f(\xi)) = 1
$$
It follows from the Chain Rule for Moebius metrics that
\begin{align*}
\frac{d f_* \rho_x}{d\rho_z}(f(\xi)) & = \frac{d f_* \rho_x}{d\rho_p}(f(\xi)) \cdot \frac{d\rho_p}{d\rho_z}(f(\xi)) \\
                                     & = \exp(B(z, p, f(\xi))) \\
                                     & = \exp(B(z, \pi(\phi_f(\gamma)), \phi_f(\gamma)(+\infty))) \\
\end{align*}

Moreover, for any $\gamma \in T^1_x X$, the same argument shows that if $\xi = \gamma(+\infty)$, then
$$
\exp(B(z, \pi(\phi_f(\gamma)), \phi_f(\gamma)(+\infty))) = \frac{d f_* \rho_x}{d\rho_z}(f(\xi))
$$
Thus
$$
\sup_{\xi \in \partial X} \frac{d f_* \rho_x}{d\rho_z}(f(\xi)) = \sup_{\gamma \in \phi_f(T^1_x X)} \exp(B(z, \pi(\gamma), \gamma(+\infty)))
$$
which gives, using the definition of the metric $d_{\mathcal M}$,
$$
\exp(d_{\mathcal M}(\rho, \rho_z)) = u_K(z)
$$
where $K = \phi_f(T^1_x X)$. Since the unique minimizer of $u_K$ is given by $y = \hat{f}(x)$, it follows that the function
$z \mapsto d_{\mathcal M}(\rho, \rho_z)$ also has a unique minimizer given by $\hat{f}(x)$. $\diamond$

\medskip

The circumcenter extension has the following naturality properties with respect to composition with isometries:

\medskip

\begin{prop} \label{natural} Let $f : \partial X \to \partial Y$ be a Moebius homeomorphism.

\smallskip

\noindent (1) If $f$ is the boundary map of an isometry $F : X \to Y$ then $\hat{f} = F$.

\smallskip

\noindent (2) If $G : X \to X, H : Y \to Y$ are isometries with boundary maps $g, h$, then
$$
\widehat{h \circ f \circ g} = H \circ \hat{f} \circ G
$$
\end{prop}

\medskip

\noindent{\bf Proof:} Let $x \in X$.

\medskip

\noindent (1) If $f$ is the boundary map of an isometry $F$, then $f_* \rho_x = \rho_{F(x)}$, so
the nearest point to $f_* \rho_x$ is $\rho_{F(x)}$, so by the previous proposition $\hat{f}(x) = F(x)$.

\medskip

\noindent (2) Note $f_* g_* \rho_x = f_* \rho_{G(x)}$. Let $z = \hat{f}(G(x))$, so $\rho_z$ is
the nearest point to $f_* \rho_{G(x)}$. Since $h_* : \mathcal{M}(\partial Y) \to \mathcal{M}(\partial Y)$ is
an isometry which preserves the subspace of visual metrics, $h_* \rho_z = \rho_{H(z)}$ is the nearest point to
$h_* f_* \rho_{G(x)} = (h \circ f \circ g)_* \rho_x$, hence by the previous proposition $H(z) = \widehat{h \circ f \circ g}(x)$, and
$H(z) = H(\hat{f}(G(x)))$ so we are done. $\diamond$

\medskip

The key to Theorem \ref{mainthm} is the following proposition:

\medskip

\begin{prop} \label{pibytwo} Let $X$ be a proper, geodesically complete CAT(-1) space. Given $\rho \in \mathcal{M}(\partial X)$,
if $x \in X$ minimizes $z \in X \mapsto d_{\mathcal M}(\rho, \rho_z)$, then for any $y \in X \cup \partial X$ distinct
from $x$, there exists $\eta \in \partial X$
maximizing $\zeta \in \partial X \mapsto \frac{d\rho}{d\rho_x}(\zeta)$ such that $\angle^{(-1)}y x \eta \geq \pi/2$.
\end{prop}

\medskip

\noindent{\bf Proof:} Let $K \subset \partial X$ be the set where $\frac{d\rho}{d\rho_x}$ attains its maximum value $e^M$,
where $M = d_{\mathcal M}(\rho, \rho_x)$, and suppose there is a
$y \in X \cup \partial X$ such that $\angle^{(-1)}y x \eta < \pi/2$ for all $\eta \in K$. Then we can choose $\epsilon, \delta > 0$ and a neighbourhood
$N$ of $K$ such that $\angle^{(-1)}y x \eta \leq \pi/2 - \epsilon$ for all $\eta \in N$, and such that $\log \frac{d\rho}{d\rho_x} \leq M - \delta$
on $\partial X - N$.

\medskip

Let $z$ be the point on the geodesic ray $[x,y)$ at a distance $t > 0$ from $x$. As $t \to 0$, for $\eta \in N$ we have,
noting that $\angle^{(-1)}z x \eta \leq \angle^{(-1)}y x \eta$,
by Lemma \ref{embedderiv},
\begin{align*}
\log \frac{d\rho}{d\rho_z}(\eta) & = \log \frac{d\rho}{d\rho_x}(\eta) - \log \frac{d\rho_z}{d\rho_x}(\eta) \\
                           & \leq M - t \cos(\angle^{(-1)}z x \eta) + o(t) \\
                           & \leq M - t \cos(\angle^{(-1)}y x \eta) + o(t) \\
                           & \leq M - t \sin \epsilon + o(t) \\
                           & < M \\
\end{align*}
for $t$ small enough depending only on $\epsilon$, while for $\eta \in \partial X - N$ we have
\begin{align*}
\log \frac{d\rho}{d\rho_z}(\eta) & = \log \frac{d\rho}{d\rho_x}(\eta) - \log \frac{d\rho_z}{d\rho_x}(\eta) \\
                           & \leq (M - \delta) + t \\
                           & < M \\
\end{align*}
 for $t < \delta$, thus for $t > 0$ small enough we have $d_{\mathcal M}(\rho, \rho_z) < M = d_{\mathcal M}(\rho, \rho_x)$,
 a contradiction. $\diamond$

\medskip

Theorem \ref{thmone} now follows from the following proposition:

\medskip

\begin{prop} \label{holder} Let $f : \partial X \to \partial Y$ be a Moebius homeomorphism
between boundaries of proper, geodesically complete CAT(-1) spaces $X, Y$. Then the circumcenter
extension $\hat{f} : X \to Y$ satisfies
$$
\cosh(d(\hat{f}(x), \hat{f}(y))) \leq e^{d(x,y)}
$$
for all $x,y \in X$. In particular $\hat{f}$ is locally $1/2$-Holder continuous: 
$$
d(\hat{f}(x), \hat{f}(y)) \leq 2 d(x,y)^{1/2}
$$
for all $x, y \in X$ such that $d(x,y) \leq 1$.
\end{prop}

\medskip

\noindent{\bf Proof:} Given $x,y \in X$, let $x' = \hat{f}(x), y' = \hat{f}(y)$. We may assume $x' \neq y'$
(otherwise the above inequality holds trivially), and also
(interchanging $x, y$ if necessary) that
$$
d_{\mathcal M}(f_* \rho_x, \rho_{x'}) \geq d_{\mathcal M}(f_* \rho_y, \rho_{y'}).
$$
Let $\rho = f_* \rho_x \in \mathcal{M}(\partial Y)$.
By Proposition \ref{nearuniq}, $x'$ minimizes $z \in Y \mapsto d_{\mathcal M}(\rho, \rho_z)$. Hence by
the previous Proposition \ref{pibytwo}, there exists $\eta \in \partial Y$
maximizing $\zeta \in \partial Y \mapsto \frac{d\rho}{d\rho_{x'}}(\zeta)$ such that $\angle^{(-1)}y' x' \eta \geq \pi/2$.
By Lemma \ref{buseform}, we have
$$
e^{B(y',x', \eta)} = \cosh(d(x', y')) - \sinh(d(x',y')) \cos(\angle^{(-1)}y' x' \eta) \geq \cosh(d(x', y'))
$$
Also,
\begin{align*}
e^{B(y',x', \eta)} & = \frac{d\rho_{x'}}{d\rho_{y'}}(\eta) \\
                   & = \frac{d\rho_{x'}}{d f_* \rho_x}(\eta) \frac{d f_* \rho_x}{d f_* \rho_y}(\eta) \frac{d f_* \rho_y}{d\rho_{y'}}(\eta) \\
                   & \leq \exp(-d_{\mathcal M}(f_* \rho_x, \rho_{x'})) \frac{d\rho_x}{d\rho_y}(f^{-1}(\eta)) \exp(d_{\mathcal M}(f_* \rho_y, \rho_{y'})) \\
                   & \leq \frac{d\rho_x}{d\rho_y}(f^{-1}(\eta)) \\
                   & = e^{B(y, x, f^{-1}(\eta))} \\
                   & \leq e^{d(x,y)} \\
\end{align*}
thus
$$
\cosh(d(x', y')) \leq e^{d(x,y)}
$$
as required.

\medskip

It follows easily that $\hat{f}$ is locally $1/2$-Holder:

\medskip

Since $e^t \leq 1 + 2t$ for $0 \leq t \leq 1$, for $x,y \in X$, if $d(x,y) \leq 1$ we have
\begin{align*}
1 + \frac{d(\hat{f}(x), \hat{f}(y))^2}{2} & \leq \cosh(d(\hat{f}(x), \hat{f}(y))) \\
                                          & \leq e^{d(x,y)} \\
                                          & \leq 1 + 2d(x,y) \\
\end{align*}
hence
$$
d(\hat{f}(x), \hat{f}(y)) \leq 2 d(x,y)^{1/2}.
$$
$\diamond$

\medskip

Let $X$ be a complete, simply connected Riemannian manifold with sectional curvatures $K$ satisfying $-b^2 \leq K \leq -1$ for
some $b \geq 1$. For $x \in X$ and $\xi, \eta \in \partial X$, let $\angle \xi x \eta \in [0, \pi]$ denote the Riemannian angle at $x$ between
the geodesic rays $[x,\xi)$ and $[x, \eta)$.

\medskip

\begin{lemma} \label{anglecomp} We have
$$
\rho_x(\xi, \eta)^b \leq \sin \left( \frac{\angle \xi x \eta}{2} \right) \leq \rho_x(\xi, \eta)
$$
\end{lemma}

\medskip

\noindent{\bf Proof:} Since the sectional curvature of $X$ is bounded above and below by $-1$ and $-b^2$, we have
$$
\angle^{(-b^2)} \xi x \eta \leq \angle \xi x \eta \leq \angle^{(-1)} \xi x \eta
$$
hence by Lemma \ref{compexist}
$$
\rho_x(\xi, \eta)^b = \sin \left( \frac{\angle^{(-b^2)} \xi x \eta}{2} \right) \leq \sin \left( \frac{\angle \xi x \eta}{2} \right)
 \leq \sin \left( \frac{\angle^{(-1)} \xi x \eta}{2} \right) = \rho_x(\xi, \eta)
$$
$\diamond$

\medskip

\noindent{\bf Proof of Theorem \ref{mainthm}}: Let $f : \partial X \to \partial Y$ be a Moebius homeomorphism
between boundaries of complete, simply connected manifolds with sectional curvatures $K$ satisfying $-b^2 \leq K \leq -1$.

\medskip

Let $x \in X$, and let $y = \hat{f}(x)$. Let $M = d_{\mathcal M}(f_* \rho_x, \rho_y)$. Let $K \subset \partial Y$ be the set
where $\frac{d f_* \rho_x}{d\rho_y}$ attains its maximum value $e^M$, and let $\eta_1 \in K$.
Then by Proposition \ref{pibytwo}, there exists
$\eta_2 \in K$ such that $\angle^{(-1)} \eta_1 y \eta_2 \geq \pi/2$, so
$\rho_y(\eta_1, \eta_2) \geq 1/\sqrt{2}$.

\medskip

Let $\xi_i = f^{-1}(\eta_i) \in \partial X, i = 1,2$.
Let $\eta'_i \in \partial Y$ be the unique point such that $\rho_y(\eta_i, \eta'_i) = 1$, $i = 1,2$. Then by Lemma
\ref{maxmin}, $\frac{d f_* \rho_x}{d\rho_y}$ attains its minimum value $e^{-M}$ at $\eta'_1, \eta'_2$, and the points
$\xi'_i = f^{-1}(\eta'_i)$ satisfy $\rho_x(\xi_i, \xi'_i) = 1$, $i = 1,2$. The Geometric Mean Value Theorem gives
$$
\rho_x(\xi_1, \xi_2) = e^M \rho_y(\eta_1, \eta_2), \rho_x(\xi'_1, \xi'_2) = e^{-M} \rho_y(\eta'_1, \eta'_2)
$$
Noting that $\angle \xi_1 x \xi_2 = \angle \xi'_1 x \xi'_2$ and $\angle \eta_1 y \eta_2 = \angle \eta'_1 y \eta'_2$, by
Lemma \ref{anglecomp} we have
\begin{align*}
\rho_x(\xi'_1, \xi'_2) & \geq \sin \left( \frac{ \angle \xi'_1 x \xi'_2}{2} \right) \\
                       & = \sin \left( \frac{ \angle \xi_1 x \xi_2}{2} \right) \\
                       & \geq \rho_x(\xi_1, \xi_2)^b \\
\end{align*}
 and
\begin{align*}
 \rho_y(\eta'_1, \eta'_1) & \leq \left(\sin \left( \frac{ \angle \eta'_1 y \eta'_2}{2} \right)\right)^{1/b} \\
                          & = \left(\sin \left( \frac{ \angle \eta_1 y \eta_2}{2} \right)\right)^{1/b} \\
                          & \leq \rho_y(\eta_1, \eta_2)^{1/b} \\
 \end{align*}
 Using the above two inequalities in the equality
 $$
 \frac{\rho_x(\xi_1, \xi_2)}{\rho_x(\xi'_1, \xi'_2)} = e^{2M} \frac{\rho_y(\eta_1, \eta_2)}{\rho_y(\eta'_1, \eta'_2)}
 $$
 gives
 $$
 \frac{1}{\rho_x(\xi_1, \xi_2)^{b-1}} \geq e^{2M} \rho_y(\eta_1, \eta_2)^{1 - 1/b}
 $$
 Thus
 \begin{align*}
 1 & \geq e^{2M} \rho_x(\xi_1, \xi_2)^{b-1} \rho_y(\eta_1, \eta_2)^{1 - 1/b} \\
  & = e^{2M} e^{(b-1)M} \rho_y(\eta_1, \eta_2)^{(b-1)+(1 - 1/b)} \\
  & \geq \frac{e^{(b+1)M}}{ \sqrt{2}^{b - 1/b} } \\
  \end{align*}
  hence
  $$
  M \leq \frac{1}{2} \frac{b - 1/b}{b+1} \log 2 = \frac{1}{2}(1 - 1/b)\log 2
  $$
  Thus
  $$
  d_{\mathcal M}(f_* \rho_x, \rho_{\hat{f}(x)}) \leq \frac{1}{2}(1 - 1/b)\log 2
  $$
  for all $x \in X$. Then for any $x, y \in X$,
  \begin{align*}
  |d(\hat{f}(x), \hat{f}(y)) - d(x, y)| & = |d_{\mathcal M}(\rho_{\hat{f}(x)}, \rho_{\hat{f}(y)}) - d_{\mathcal M}(f_* \rho_x, f_* \rho_y)| \\
                                        & \leq d_{\mathcal M}(f_* \rho_x, \rho_{\hat{f}(x)}) + d_{\mathcal M}(f_* \rho_y, \rho_{\hat{f}(y)}) \\
                                        & \leq (1 - 1/b) \log 2 \\
  \end{align*}
  thus $\hat{f}$ is a $(1, (1 - 1/b)\log 2)$-quasi-isometry. As in \cite{biswas3} it is straightforward to show that the
  image of $\hat{f}$ is $\frac{1}{2}(1 - 1/b)\log 2$-dense in $Y$ and that the boundary map of $\hat{f}$ equals $f$. $\diamond$

\bibliography{moeb}
\bibliographystyle{alpha}

\end{document}